\documentclass[12pt]{article}
\usepackage{amssymb,amsmath, amsthm, amscd, ifthen}

\usepackage{graphicx,url}

\usepackage[left=2cm, top=2cm, right=2cm, bottom=20mm, nohead]{geometry}

\usepackage{authblk}

\newsavebox\affbox

\usepackage[russian]{babel}   
\usepackage[utf8]{inputenc}  

\newcommand{\Natur}{{\mathbb N}}

\newtheorem{lemma}{\emph{Лемма}}

\newtheorem{theorem}{Теорема}
\newtheorem{claim}{Утверждение}

\title{Оценка числа рёбер в особых подграфах некоторого дистанционного графа \thanks{Настоящая работа выполнена за счет гранта РФФИ (проект N 18-01-00355) и гранта президента НШ-6760.2018.1.}}

\author[1]{Ф.А. Пушняков\thanks{filipp.pushnyakov@phystech.edu}}
\author[1, 2]{А.М. Райгородский\thanks{mraigor@yandex.ru}}
\affil[1]{Московский физико-технический институт (государственный университет), кафедра
дискретной математики и лаборатория продвинутой комбинаторики и сетевых приложений}
\affil[2]{МГУ им. М.В. Ломоносова, механико-математический факультет, кафедра математической статистики и случайных 
процессов; Адыгейский государственный университет, Кавказский математический центр; Бурятский государственный университет, институт математики и информатики.
}

\begin{document} 

\maketitle

\section{Введение}

Рассмотрим последовательность графов $G_{n} = G_{n} (V_{n}, E_{n}) = G(n,3,1)$,  у которых $$V_{n} = \{x = (x_{1}, \dots, x_{n}) \; \vert \; x_{i} \in \{0, 1\},\; i = 1, \dots , n \; , \; x_{1} + \ldots + x_{n} = 3\},$$ $$E_{n} = \{ (x, y) \; \vert \; \langle x, y \rangle = 1\},$$
где через $\langle x, y \rangle$ обозначено скалярное произведение векторов $x$ и $y$. 
Иными словами, вершинами графа $G\left(n, 3, 1\right)$ являются $\left(0,1\right)$-векторы, скалярный квадрат которых равен трем. И эти вершины соединены ребром тогда и только тогда, когда скалярное произведение соответствующих векторов равно единице. Данное определение можно переформулировать в комбинаторных терминах. А именно, рассмотрим граф, вершинами которого являются всевозможные трехэлементные подмножества множества $\mathcal{R}_{n} = \{1, \dots, n\}$, причем ребро между такими вершинами проводится тогда и только тогда, когда соответствующие трехэлементные подмножества имеют ровно один общий элемент. Изучение данного графа обусловлено многими задачами комбинаторной геометрии, экстремальной комбинаторики, теории кодирования: например, задачей Нелсона--Эрдёша--Хадвигера о раскраске метрического пространства (см. \cite{Rai3}--\cite{Rai15}), проблемой Борсука о разбиении пространства на части меньшего диаметра (см. \cite{Rai3}--\cite{Rai5}, \cite{Bolt}--\cite{Rai2}), задачами о числах Рамсея (см. \cite{Ram}--\cite{Rai12}), задачами о кодах с одним запрещенным расстоянием (см. \cite{MS}--\cite{Bobu3}).
\par Напомним несколько свойств данного графа. Граф $G(n,3,1)$ является регулярным со степенью вершины $d_{n} = 3 \cdot C_{n-3}^{2}$. Очевидно, что $\vert V_{n} \vert = C_{n}^{3} \sim \frac{n^{3}}{6}$ при  $n \rightarrow \infty$. В силу регулярности рассматриваемого графа имеем $\vert E_{n} \vert = \frac{d_{n} \cdot \vert V_{n} \vert}{2} = \frac{3}{2} \cdot C_{n-3}^{2} \cdot C_{n}^{3} \sim \frac{n^{5}}{8}$ при $n \rightarrow \infty$. \par
Напомним, что \textit{независимым множеством} графа называется такое подмножество его вершин, что никакие две вершины подмножества не соединены ребром. \textit{Числом независимости} $\alpha(G)$ называется наибольшая мощность независимого множества.  Положим $\alpha_{n} = \alpha(G(n, 3, 1))$. Результат теоремы Ж. Надя (см. \cite{Nagy}) отвечает на вопрос о числе независимости графа $G(n, 3, 1)$. А именно, $\alpha_{n} \sim n$ при $n \rightarrow \infty$. Более того, из доказательства теоремы Ж. Надя можно сделать вывод о структуре независимого множества в рассматриваемом графе. Для описания этой структуры введем дополнительные обозначения.  Пусть $W \subseteq V_{n}$. Будем говорить, что $W$ является множеством вершин \textit{первого типа}, если $\vert W \vert \geq 3$ и существуют такие $i, j \in \mathcal{R}_{n}$, что для любой вершины $w \in W$ выполнено $i, j \in w$; далее, $W$ является множеством вершин \textit{второго типа}, если $\vert W \vert \geq 2$ и существуют такие $i, \; j, \; k, \; t \in \mathcal{R}_{n}$, что для любой вершины $w \in W$ выполнено $w \subset \{i, j, k, t\} $; наконец, $W$ является множеством вершин \textit{третьего типа}, если для любых $w_{1}, w_{2} \in W$ выполнено соотношение $w_{1} \cap w_{2} = \emptyset$. Более того, \textit{носителем} множества вершин назовем объединение всех вершин данного множества. Тогда имеет место следующее утверждение.
\begin{claim} Любое независимое множество $U \subseteq V_{n}$ можно представить в виде объединения $$U = \left( \cup_{i \in \mathcal{I}} A_{i} \right) \cup \left( \cup_{j \in \mathcal{J}} B_{j} \right) \cup \left( \cup_{k \in \mathcal{K}} C_{k} \right),$$ где $A_{i}$ -- множество вершин первого типа, $B_{j}$ -- множество вершин второго типа, $C_{k}$ -- множество вершин третьего типа, $i \in \mathcal{I}, \; j \in \mathcal{J},\;  k \in \mathcal{K}$, и носители всех упомянутых множеств попарно не пересекаются.
\end{claim}
Мы не доказываем данное утверждение, так как оно мгновенно следует из доказательства теоремы Ж. Надя (см. \cite{Nagy}).
\par
Обозначим через $r(W)$ количество рёбер графа $G$ на множестве $W \subseteq V_{n}$. Иными словами, $$r(W) = \vert \{(x, y) \in E(G) \; \vert \; x \in W, \; y \in W\} \vert \; .$$ Также положим  $$r(l(n)) = \min_{\vert W \vert = l(n), \; W \subseteq V_{n}} r(W) \; .$$
Заметим, что если $l(n) \leq \alpha_{n}$, то $r(l(n)) = 0$ и обсуждать нечего. Если же $l(n) > \alpha_{n}$, то, очевидно, в любом $W \subseteq V_{n}$ мощности $l(n)$ непременно найдутся рёбра. Возникает интересный вопрос об изучении величины $r(l(n))$. Оценки, полученные в работах \cite{Rai11}--\cite{Rai12}, достаточно слабые, поэтому появились работы \cite{Pushnyakov}--\cite{Pushnyakov2}, в которых приведено практически полное исследование величины $r(l(n))$. А именно, в работе \cite{Pushnyakov} была доказана следующая теорема (читая формулировку, важно помнить, что $ n \sim \alpha_n $).

\begin{theorem} \label{t1}
Имеют место четыре случая:
\begin{enumerate}
\item Пусть функции $f: \Natur \rightarrow \Natur, \; g: \Natur \rightarrow \Natur$ таковы, что выполнено $n = o(f)$ и $g = o(n^{2})$ при $n \rightarrow \infty$. Пусть функция $l: \Natur \rightarrow \Natur$ такова, что для любого $n \in \Natur$  выполнена цепочка неравенств $f(n) \leq l(n) \leq g(n)$. Тогда $r(l(n)) \sim \frac{l(n)^{2}}{2 \alpha_{n}}$ при $n \rightarrow \infty$.
\item Пусть функция $l: \Natur \rightarrow \Natur$  такова, что существуют константы $C_{1}, \; C_{2}$, с которыми для каждого $n \in \Natur$ выполнена цепочка неравенств $C_{1} \cdot n^{2} \leq f(n) \leq C_{2} \cdot n^{2}$. Тогда $r(l(n)) \sim \frac{l(n)^{2}}{2 \alpha_{n}}$ при $n \rightarrow \infty$.
\item Пусть функции $f: \Natur \rightarrow \Natur, \; g: \Natur \rightarrow \Natur$ таковы, что выполнено $n^{2} = o(f)$ и $g = o(n^{3})$ при $n \rightarrow \infty$.  Пусть функция $l: \Natur \rightarrow \Natur$ такова, что для каждого $n \in \Natur$ выполнено $f(n) \leq l(n) \leq g(n)$. Тогда существуют такие функции $h_{1}: \Natur \rightarrow \Natur, \ h_{2}: \Natur \rightarrow \Natur$, что $h_{1} \sim \frac{l(n)^{2}}{\alpha_{n}}, \ h_{2} \sim \frac{5 l(n)^{2}}{\alpha_{n}}$  при $n \rightarrow \infty$ и для каждого $n \in \Natur$ выполнена цепочка неравенств $h_{1}(n) \leq r(l(n)) \leq h_{2}(n)$. При этом для выполнения нижней оценки требование $ g = o(n^3) $ не нужно.
\item Пусть функция $l: \Natur \rightarrow \Natur$  такова, что существует константа $C$, с которой выполнена цепочка неравенств $C \cdot n^{3} \leq l(n) \leq  C_{n}^{3}$. Пусть $c_{n} = 1-\frac{l(n)}{C_{n}^{3}}$. Тогда существует функция $f: \Natur \rightarrow \Natur$, такая, что $f(n) \sim n^{5} \left(\frac{1}{8} - \frac{c_{n}}{4} + \frac{c_{n}^{2}}{72}\right)$ при $n \rightarrow \infty$ и для каждого $n \in \Natur$ выполнено $r(l(n)) \geq f(n)$.
\end{enumerate}
\end{theorem}

Как можно заметить, оценки, полученные в пунктах 3-4 данной теоремы, не являются точными. Отметим в то же время, что нижняя оценка из пункта 3 верна и в условиях пункта 4. В работе \cite{Pushnyakov2} была улучшена оценка из пункта 3. А именно, была доказана следующая теорема.

\begin{theorem} \label{t3}
Пусть функция $l: \Natur \rightarrow \Natur$ такова, что $n^{2} = o(l)$ при $n \rightarrow \infty$. Тогда существует такая функция $h: \Natur \rightarrow \Natur$, что $h \sim \frac{3 l^{2}}{2n}$ при $n \rightarrow \infty$ и $r(l(n)) \geq h(n)$ для достаточно большого $n$. 
\end{theorem}

Более того, в работе \cite{Pushnyakov3} первый автор устранил одну неточность в пункте 4 теоремы 1 и улучшил многие оценки из теоремы 1, а именно доказал следующую теорему.

\begin{theorem} \label{t2}
Имеют место четыре случая:
\begin{enumerate}
\item
Пусть дана произвольная функция $l: \Natur \rightarrow \Natur$ с ограничением $n = o(l)$. Тогда существует такая функция $h: \Natur \rightarrow \Natur$, что $h(n) \sim \frac{9 l(n)^{2}}{2 \alpha_{n}} $  при $n \rightarrow \infty$ и для каждого $n \in \Natur$ выполнено неравенство $r(l(n)) \leq h(n)$.
\item Пусть функция $l: \Natur \rightarrow \Natur$  такова, что существуют константа $C$ и функция $g: \Natur \rightarrow \Natur$ такая, что $n^{2} = o(g)$ при $n \rightarrow \infty$ и для каждого $n \in \Natur$ выполнена цепочка неравенств $C \cdot n^{3} \leq l(n) \leq  C_{n}^{3} - g(n)$. Пусть $c_{n} = 1-\frac{l(n)}{C_{n}^{3}}.$ 
Тогда существует такая функция $h : \Natur \rightarrow \Natur$, что $h = o\left(1\right)$ при $n \rightarrow \infty$ и для любого достаточно большого $n \in \Natur$ выполнено неравенство $r(l(n)) \geq \frac{n^{5}}{8} \left(1 - 2c_{n} + \frac{c_{n}^{2}}{3} \left(1 + h\left(n\right)\right) - \frac{10}{n} + \frac{20 \cdot c_{n}}{n} - \frac{10 \cdot c_{n}^{2}}{3 n} \left(1 + h\left(n\right)\right)\right).$
\item Пусть функция $l: \Natur \rightarrow \Natur$  такова, что существуют константы $B, \ C$ и функция $g: \Natur \rightarrow \Natur$ такая, что $n = o\left(g\right)$ при $n \rightarrow \infty$ и для каждого $n \in \Natur$ выполнено $g(n) \leq Bn^{2}$ и $C \cdot n^{3} \leq l(n) \leq  C_{n}^{3} - g(n)$. Пусть $c_{n} = 1-\frac{l(n)}{C_{n}^{3}}.$ 
Тогда существует такая функция $h : \Natur \rightarrow \Natur,$ что $h = o\left(1\right)$ при $n \rightarrow \infty$ и для любого достаточно большого $n \in \Natur$ выполнено неравенство $r(l(n)) \geq \frac{n^{5}}{8} \left(1 - 2 c_{n} + \frac{2 c_{n}^{2}}{9} \left(1 + h\left(n\right)\right) - \frac{10}{n} + \frac{20 \cdot c_{n}}{n} - \frac{20 \cdot c_{n}^{2}}{9 n} \left(1 + h\left(n\right)\right)\right).$
\item Пусть функция $l: \Natur \rightarrow \Natur$  такова, что существует такая константа $C$, что для каждого $n \in \Natur$ выполнена цепочка неравенств $C_{n}^{3} - Cn \leq l(n) \leq  C_{n}^{3}$. Пусть $c_{n} = 1-\frac{l(n)}{C_{n}^{3}}$. 
Тогда для любого $n \in \Natur$ выполнено неравенство $r(l(n)) \geq \frac{n^{5}}{8} \left(1 - 2 c_{n} - \frac{10}{n} + \frac{20 \cdot c_{n}}{n}\right).$

\end{enumerate}
\end{theorem}
Результат первого пункта теоремы \ref{t2} улучшает оценку, полученную в пункте 3 теоремы \ref{t1}. Тем не менее, эта оценка по-прежнему не является точной: величина $r(l(n))$ удовлетворяет следующей цепочке неравенств: 
$$\frac{3 l(n)^{2}}{2 \alpha_{n}} (1 + o(1)) \leq r(l(n)) \leq \frac{9 l(n)^{2}}{2 \alpha_{n}} (1 + o(1))$$
при $n \rightarrow \infty$. Между нижней и верхней оценками имеется зазор в 3 раза.
\par 
Результат второго и третьего пунктов теоремы \ref{t2} немного улучшает оценку, полученную в пункте 4 теоремы \ref{t1}. 
\par Посмотрим ещё с несколько иной точки зрения на полученные результаты. Можно записать лучшую известную нам верхнюю оценку в виде 
\begin{equation}\label{formula1}
	r(l(n)) \le \frac{9l(n)^{2}}{2n} (1+o(1)) = \frac{n^{5}}{8} (1-2c_{n}+c_{n}^{2}) (1+o(1)),
\end{equation}
где $ c_n $ из формулировки теоремы 3. В таком же виде можно записать и нижнюю оценку из теоремы 2:
\begin{equation}\label{formula2}
	r(l(n)) \ge \frac{n^5}{8} \left(\frac{1}{3}-\frac{2}{3}c_{n}+\frac{1}{3} c_{n}^2\right) (1+o(1)).
\end{equation}

Как мы помним, в условиях пункта 4 теоремы 1 она верна, то есть верна она и в условиях пунктов 2--4 теоремы 3. Конечно, если $ c_{n} \rightarrow 0 $, то оценки пунктов 2--4 новой теоремы асимптотически совпадают с оценкой (\ref{formula1}) и в этом случае оценка (\ref{formula2}) им не конкурент. Однако в условиях теоремы 3 возможно и что $ c_{n} $ не стремится к нулю (хотя и не превосходит константы, строго меньшей единицы). В этом случае оценки из пунктов 2--4 становятся лучше, чем оценка из теоремы 2 при выполнении неравенства  
$$\frac{n^{5}}{8} \left(1 - 2 c_{n} - \frac{10}{n} + \frac{20 \cdot c_{n}}{n}\right) \geq \frac{3 \left( c_{n} C_{n}^{3}\right)^{2}}{2},$$
которое выполнено при $c_{n} \le 0.486 \dots$ и достаточно больших значениях $n \in \Natur$. 

\par 
Заметим, что все оценки, приведённые выше, выполнены для всех возможных подграфов графа $G\left(n, 3, 1\right)$. Возникает вопрос, а можно ли улучшить оценки, если рассмотреть только лишь подграфы определённой структуры? Например, подграфы, не содержащие какой-нибудь заданной конструкции. В данной работе автору удалось значительно улучшить оценки величины $r\left(l\left(n\right)\right)$ для некоторого класса подграфов графа $G(n, 3, 1)$. Для формулировки этих результатов нам потребуется определить некоторые дополнительные понятия. 

\textit{Звездным множеством} будем называть любое независимое множество, состоящее из вершин, принадлежащих множествам первого и третьего типов. 
Диаметром $d\left(A\right)$ звездного множества $A$ назовём мощность носителя $A$. И, наконец, диаметром $d\left(W\right)$ множества вершин $W$ назовём максимальный диаметр звездного множества, содержащегося в $W$. Формально говоря, 
$$d\left(W\right) = \max_{A \subset W} \{ d\left(A\right) \; \vert \; A \text{ является звездным множеством}\}.$$

Также для произвольной функции $\rho : \Natur \rightarrow \Natur$, удовлетворяющей неравенству $\rho\left(n\right) \leq n$ для любого натурального $n$, и произвольного $W \subset V_{n}$ положим 
\begin{equation*}
r_{\rho}\left(W\right) = \begin{cases}
r\left(W\right) &\text{\text{при } \; $d\left(W\right) \leq \rho\left(n\right)$}\\
C_{n}^{5} &\text{\text{при } \; $d\left(W\right) > \rho\left(n\right)$}\\
\end{cases}
\end{equation*}

Рассмотрим подробнее определение функции $r_{\rho}$. Легко заметить, что для всех $W \subset V_{n}$ выполнено неравенство $r_{\rho}\left(W\right) \geq r\left(W\right)$, а любые нижние оценки являются либо тривиальными, либо оценками для подмножеств вершин $W$, у которых диаметр не превосходит некоторой функции $\rho$. 
\par 
Основным результатом данной статьи является следующая теорема.

\begin{theorem}
	Пусть функции $l : \Natur \rightarrow \Natur, \; \rho : \Natur \rightarrow \Natur$ таковы,  что $n^{2} = o\left(l\right)$ при $n \rightarrow \infty$, и для любого натурального $n$ выполнено неравенство $l(n) \leq C_{n}^{3}$. Тогда для любого $W \subset V_{n}$ мощности $l\left(n\right)$ выполнено неравенство
	\begin{equation}\label{formula3}
	r_{\rho}\left(W\right) \geq \frac{l^{2}}{n} \left(2 - \frac{\rho\left(n\right)^{3}}{6l} + o\left(1\right)\right)
	\end{equation}
	при $n \rightarrow \infty$.  
\end{theorem}

\par
Как уже было замечено, теорема 4 не всегда является улучшением старых результатов. Более того, в некоторых случаях оценка из теоремы 4 является тривиальной. А именно, в случае $d\left(W\right) > \rho\left(n\right)$ доказывать нечего --- число рёбер в произвольном подмножестве вершин графа $G\left(n, 3, 1\right)$ никак не может быть больше общего числа рёбер графа. С другой стороны, при $d\left(W\right) \leq \rho\left(n\right)$ оценка является нетривиальной, и её-то мы и будем доказывать! 
\par 
Ясно, что при $\rho^3 = o\left(l\right)$ результат теоремы 4 является значимым улучшением по сравнению со всеми прежними результатами, ведь константа $\frac{3}{2}$ в них заменена константой 2. При $l$ порядка $n^{3}$ данное условие состоит в малости $\rho$ в сравнении с $n$. Но даже если $\rho$ порядка $n$, новая оценка может быть сильнее прежних. Например, если $l = \frac{C_{n}^{3}}{2}$, а $\rho = \frac{n}{2}$, то правая часть оценки из теоремы может быть записана в виде $\frac{l^{2}}{n} \left(1.75 - \dots \right)$, что является несомненным улучшением старых результатов.
 Тем не менее, оценка из теоремы 4 по-прежнему далека от лучшей верхней оценки. Действительно, известные нам результаты можно упрощённо записать в виде 
$$ \frac{9 l^{2}}{2 \alpha_{n}} \left( 1 + o\left(1\right)\right) \geq r\left(W\right) \geq \frac{l^{2}}{n} \cdot \left(2 - \dots \right).$$
Как легко заметить, зазор между левой и правой частями неравенства по-прежнему существенный.

\section{Доказательство теоремы 4} \label{sec:firstpage}
\subsection{Вспомогательные утверждения и определения}
Перед началом доказательства введём вспомогательное определение и сформулируем несколько вспомогательных утверждений. Для произвольного множества вершин $S$ и произвольной вершины $v \notin S$ положим 
$$n\left(v, S\right) = \left| \{ u \in S \; \vert \; \left(u, v\right) \in E_{n}\} \right|.$$
Иными словами, $n(v, S)$ обозначает число вершин множества $S$, соединённых ребром с вершиной $v$. 
\par Пусть $H$ -- произвольное подмножество вершин графа $G\left(n, 3, 1\right)$. Пусть также $I$ -- наибольшее независимое подмножество подграфа графа $G\left(n, 3, 1\right)$, индуцированного множеством вершин $H$. Очевидно, $\vert I \vert \leq  n$. Положим 
$$B_{i} = \{ w \in H \ \vert \ n\left( w, I\right) = i\}.$$
Иными словами, $B_{i}$ -- это подмножество вершин графа $H$, которые соединены ровно с $i$ вершинами из множества $I$. Оценим мощности множеств $B_{i}$ для некоторых значений $i$. 
\par 
Очевидно, что $\vert B_{0}\vert = 0$, так как иначе $I$ не было бы максимальным независимым множеством. 
В статье \cite{Pushnyakov2} была доказана следующая лемма.
\begin{lemma}
	В обозначениях выше выполнено неравенство $\vert B_{1} \vert + \vert B_{2} \vert \leq 35 \cdot n^{2}$.
\end{lemma}
Таким образом, число вершин, соединённых с не более, чем двумя вершинами независимого множества $I$, достаточно маленькое. Оказывается, число вершин, соединённых ровно с тремя вершинами множества $I$, также не является достаточно большим. А именно, верна следующая лемма.
\begin{lemma}
	$\left| B_{3} \right| \leq \frac{\left(\rho\left(n\right)\right)^{3}}{6} + 20n^{2}$.
\end{lemma}

Доказательство леммы 2 будет приведено в пункте 4.  

\subsection{Доказательство теоремы}

\par Теперь перейдём к доказательству теоремы. Пусть $I_{1}$ -- наибольшее независимое подмножество вершин множества $W$. Положим $\alpha_{1} = \vert I_{1} \vert$. Ясно, что $\alpha_{1} \leq n$. Рассмотрим множество $W_{1} = W \setminus I_{1}$. Каждая вершина из данного множества соединена как минимум с одной вершиной из множества $I_{1}$, иначе $I_{1}$ не было бы максимальным. Положим,
$$f_{1} = \vert \{ v \in W_{1} \; \vert \; n\left(v, I_{1}\right) \leq 2\} \vert.$$
 Иными словами, $f_{1}$ -- это мощность множества вершин из множества $W_{1}$, соединенных с не более чем двумя вершинами из $I_{1}$. Как следует из леммы 1, $f_{1} \leq 35 \cdot  n^{2}$.
Также, пусть $f_{2}$ обозначает число вершин, соединенных ровно с тремя вершинами из множества $I_{1}$. Формально говоря, 
$$f_{2} = \vert \{ v \in W_{1} \; \vert \; n\left(v, I_{1}\right) = 3\} \vert.$$
По лемме 2 выполнено неравенство $f_{2} \leq \frac{\left( \rho \left(n\right)\right)^{3}}{6} + 20n^{2}$. Остальные вершины из множества $W_{1}$ соединены хотя бы с четырьмя вершинами из множества $I_{1}$. Стало быть, число рёбер между вершинами множеств $W_{1}$ и $I_{1}$ можно оценить снизу как 
$$4\left(l - \alpha_{1} - f_{1} - f_{2}\right) + 3 f_{2} + 2 f_{1} = $$
$$= 4 l - 4 \alpha_{1} - 2f_{1} - f_{2} \geq 4 l - 4 n - 90 \cdot n^{2} - \frac{\left(\rho\left(n\right)\right)^{3}}{6}.$$ \par 
Теперь попробуем повторить подобную операцию несколько раз. А именно рассмотрим множество вершин $W_{1}$. Выберем в нём максимальное независимое множество и обозначим его $I_{2}$. Аналогично положим $\alpha_{2} = \vert I_{2} \vert$. Ясно, что $\alpha_{2} \leq \alpha_{1} \leq n$, так как на первом шаге мы взяли наибольшее независимое множество. Обозначим $W_{2} = W_{1} \setminus I_{2}$. Опять же, как и в предыдущем случае, каждая из вершин множества $W_{2}$ соединена хотя бы с одной вершиной из множества $I_{2}$. Положим 
$$f_{3} = \vert \{ v \in W_{2} \; \vert \; n\left(v, I_{2}\right) \leq 2\} \vert.$$
Иными словами, $f_{3}$ -- это мощность множества вершин из $W_{2}$, соединённых с не более чем двумя вершинами из множества $I_{2}$. 
И, наконец, через $f_{4}$ обозначим число вершин из множества $W_{2}$, соединённых ровно с тремя вершинами из $I_{2}$. Формально говоря, 
$$f_{4} = \vert \{ v \in W_{2} \; \vert \; n\left(v, I_{2}\right) = 3\} \vert.$$
Из леммы 2 следует, что  $f_{4} \leq \frac{\left(\rho\left(n\right)\right)^{3}}{6} + 20n^{2}$. Остальные вершины из множества $W_{2}$, очевидно, соединены хотя бы с четырьмя вершинами из $I_{2}$. Таким образом, количество рёбер между вершинами из множеств $W_{2}$ и $I_{2}$ можно оценить снизу как 
$$4\left(l - \alpha_{1} - \alpha_{2} - f_{3} - f_{4}\right) + 3 f_{4} + 2 f_{3} = $$
$$= 4 l - 4 \alpha_{1} - 4 \alpha_{2} - 2f_{3} - f_{4} \geq 4 l - 8 n - 90 n^{2} - \frac{\left(\rho\left(n\right)\right)^{3}}{6}.$$
\par 
Продолжим данный процесс 
$$t = \left[ \frac{l}{n} \right]$$ шагов. Получим, что мы найдём как минимум
$$\sum_{i=1}^{t} \left(4l - 4in - 90 \cdot n^{2} - \frac{\left(\rho\left(n\right)\right)^{3}}{6}\right)$$
рёбер. Оценим каждое слагаемое данной суммы отдельно, а потом из этих оценок получим оценку всей суммы. 

$$\sum_{i=1}^{t} 4l = 4 l \cdot t \geq 4l \left(\frac{l}{n} - 1\right) = \frac{4l^{2}}{n} - 4l.$$

$$\sum_{i=1}^{t} 4in = 4n \frac{t\left(t + 1\right)}{2} \leq \frac{4n}{2}\frac{l}{n}\left(\frac{l}{n}+1\right) =\frac{2l^{2}}{n} + 2l.$$

$$\sum_{i=1}^{t} 90 \cdot n^{2} = 90 \cdot tn^{2} \leq 90 \cdot n^{2}\frac{l}{n} = 90 nl.$$

$$\sum_{i=1}^{t} \frac{\rho\left(n\right)^{3}}{6} = \frac{t \cdot \rho\left(n\right)^{3}}{6} \leq \frac{\rho\left(n\right)^{3}l}{6n}.$$

Просуммируем полученные неравенства и получим, что число рёбер, найденное нами в результате описанной выше процедуры, не меньше, чем 
$$ \frac{l^{2}}{n} \left(2 - \frac{\rho\left(n\right)^{3}}{6l} - 90\frac{n^{2}}{l} - 6\frac{n}{l}\right) = \frac{l^{2}}{n} \left(2 - \frac{\rho\left(n\right)^{3}}{6l} + o\left(1\right)\right)$$ при $n \rightarrow \infty$. Таким образом, утверждение теоремы доказано.

\section{Доказательство леммы 2}

\par 
Пусть $w$ -- произвольная вершина, принадлежащая множеству $B_{3}$. По определению, она соединена ровно с тремя вершинами независимого множества $I$. В силу утверждения 1 вершины независимого множества $I$ можно разбить на три непересекающихся подмножества: подмножество, все вершины которого лежат в некотором множестве вершин первого типа, подмножество вершин, все вершины которого лежат в некотором множестве второго типа, и, наконец, подмножество вершин, все вершины которого лежат в некотором множестве третьего типа. 
\par
Введём дополнительные обозначения. Пусть  
$$F = \bigcup_{t \in \mathcal{F}}\{F_{t}\} \text{ -- множество всех уникальных подмножеств вершин } I \text{ первого типа,}$$ 
$$S = \bigcup_{t \in \mathcal{S}}\{S_{t}\} \text{ -- множество всех уникальных подмножеств вершин } I \text{ второго типа,}$$
$$T = \bigcup_{t \in \mathcal{T}}\{T_{t}\} \text{ -- множество всех уникальных подмножеств вершин } I \text{ третьего типа.}$$
\par
Напомним, что носители всех упомянутых множеств попарно не пересекаются. \par 

\par 
Введём два вспомогательных определения. Напомним, что множеством вершин второго типа является множество вершин $S$, для которого существуют такие четыре различных элемента $\{i, j, k, l\} \in \mathcal{R}_{n}$, что носители вершин множества $S$ являются подмножествами множества $\{i, j, k, l\}$. Очевидно, что для фиксированных элементов $\{i, j, k, l\}$ существует ровно четыре уникальные вершины графа $G(n, 3, 1)$, носители которых лежат в множестве $\{i, j, k, l\} $. Так что совершенно ясно, что множество вершин второго типа может иметь мощность как минимум 3 и как максимум 4. Таким образом, будем называть множество вершин второго типа \textit{полным}, если его мощность равняется четырём, а иначе, если его мощность равняется трём, то \textit{неполным}.
\par 
Аналогично, будем называть элемент $a$, принадлежащий носителю множества вершин второго типа, \textit{полным}, если он принадлежит трём вершинам множества вершин второго типа, а иначе \textit{неполным}.
\par 
Перейдём к доказательству леммы. Пусть вершина $w$ пересекается с вершинами $v_{1}, v_{2}, v_{3}$ из множества $I$. Каждая из этих вершин принадлежит какому-либо подмножеству вершин либо первого, либо второго, либо третьего типа. 
Рассмотрим все возможные случаи подобной принадлежности.

\begin{enumerate}
\item \textbf{Существует такое полное множество вершин $A$ второго типа, что $\vert supp\left(A\right) \cap supp\left(w\right) \vert = 1$. }
\par
Иными словами, $w$ пересекается с каким-то полным множеством вершин по одному элементу. Пусть $x = supp\left(A\right) \cap supp\left(w\right)$. Пусть также $w = \{x, y, z\}$. В таком случае элемент $y$ можно выбрать $n-4$ способами, а оставшийся элемент $z$ не более чем четырьмя способами, так как иначе мы нашли бы независимое множество большей мощности, чем исходное. Элемент $x$ можно выбрать 4 способами, а само множество $A$ можно выбрать не более чем $\frac{n}{4}$ способами. Таким образом, в этом случае мы имеем не больше $4 \cdot \frac{n}{4} \cdot \left(n - 4\right) \cdot 4 \leq 4n^{2}$ вершин, удовлетворяющих данному условию.
\item \textbf{Существует такое полное множество вершин $A$ второго типа, что $\vert supp\left(A\right) \cap supp\left(w\right) \vert = 2$.}
\par Иными словами, существует такое полное множество вершин второго типа, с которым вершина $w$ пересекается ровно по двум элементам.
Пусть множество $A$ имеет носитель $\{x, y, z, t\}$, а $w$ имеет вид $\{x, y, a\}$. Пару элементов $\{x, y\}$ можно выбрать $C_{4}^{2}$ способами, а оставшийся элемент $a$ можно выбрать $n-4$ способами. Само множество $A$ можно выбрать не более чем $\frac{n}{4}$ способами. Таким образом, в этом случае мы имеем не больше $\frac{n}{4} \cdot C_{4}^{2} \cdot \left(n - 4\right) \leq 2n^{2}$ вершин, удовлетворяющих данному условию.
\item \textbf{Существует такое неполное множество вершин $A$ второго типа, что $\vert supp\left(A\right) \cap supp\left(w\right) \vert = 1.$} 
\par Иными словами, существует такое неполное множество вершин, с которым вершина $w$ пересекается ровно по одному элементу. Обозначим этот элемент $x$. Тогда имеют место два случая: 
\begin{itemize}
\item \textbf{Элемент $x$ является полным.}
Пусть неполное множество вершин $A$ имеет вид $\{\{x, y, z\}, \{x, y, t\}, \{x, z, t\}\}$, а вершина $w$ имеет вид $\{x, a, b\}$. Элемент $a$ можно выбрать не более чем $n-4$ способами, а элемент $b$ при фиксированном элементе $a$ можно выбрать не более чем четырьмя способами (иначе выбранное независимое множество не является максимальным). Таким образом, в данном случае мы имеем не более чем $\frac{n}{4} \cdot \left(n - 4 \right) \cdot 4 \leq n^{2}$ вершин, удовлетворяющих условию.
\item \textbf{Элемент $x$ является неполным.}	
Пусть неполное множество вершин $A$ имеет вид $\{\{x, y, z\}, \{x, y, t\}, \{y, z, t\}\}$, а вершина $w$ имеет вид $\{x, a, b\}$. В таком случае, вершина $w$ пересекается с двумя вершинами из $A$ и с одной вершиной из остального независимого множества.  
Пусть существует такая вершина $u_{1} \in I$, для которой $supp \left(u_{1}\right) \cap supp \left(w\right) = \{a\}$, и не существует такой вершины $u_{2}$, для которой $supp \left(u_{2}\right) \cap supp \left(w\right) = \{b\}$. Элемент $a$ можно выбрать не более, чем $n - 4$ способами. А элемент $b$ можно выбрать не более, чем 5 способами, так как иначе мы бы могли найти большее независимое множество, чем мы выбрали изначально. А именно, если бы существовали элементы $b_{1}, b_{2}, \dots, b_{5}$, удовлетворяющие условиям выше, то множество  
$$ \left(I \cup_{i=1}^{5} \{\{x, a, b_{i}\}\}\right) \setminus \left(\{\{x, y, z\}\} \cup \{\{x, y, t\}\} \cup \{\{y, z, t\}\} \cup \{u_{1}\}\right) $$
было бы независимым и имело бы большую мощность, чем $I$. Само множество $A$ можно выбрать не более, чем $n-4$ способами. 
Таким образом, в данном случае мы нашли не больше, чем $\frac{n}{4} \cdot \left(n-4\right) \cdot 5 \leq 2n^{2}$ вершин.
\end{itemize}

\item \textbf{Существует такое неполное множество вершин $A$ второго типа, что $\vert supp\left(A\right) \cap supp\left(w\right) \vert = 2.$ }
\par Иными словами, $w$ пересекается с неполным множеством вершин $A$ по двум элементам.
Данные 2 элемента можно выбрать не более, чем $C_{4}^{2}$ способами. Оставшийся третий элемент вершины $w$ можно выбрать не более, чем $n-4$ способами. А само множество $A$ можно выбрать не более, чем $\frac{n}{4}$ способами. Таким образом, мы имеем не более $\frac{n}{4} \cdot C_{4}^{2} \cdot \left(n - 4\right) \leq 2n^{2}$ вершин, удовлетворяющих данному условию.
\item \textbf{Не существует множества вершин $A$ второго типа, для которого $\vert supp\left(A\right) \cap supp\left(w\right) \vert > 0$.}
\par Иными словами, вершина $w$ не пересекается ни с полным, ни с неполным множеством вершин.
В данном случае $w$ пересекается только с вершинами первого и третьего типов. Возможны два случая: либо носитель вершины $w$ полностью лежит в объединении носителей вершин первого и третьего типов, либо нет. Обозначим $X_{n} = supp\left(\cup_{a \in F \cup T} \; \{a\}\right)$. Иными словами, $X_{n}$ -- это носитель объединения вершин первого и третьего типов.
\begin{itemize}
\item \textbf{$supp\left(w\right) \subset X_{n}$}
\par Иными словами, носитель вершины $w$ лежит в объединении носителей вершин первого и третьего типа. Ясно, что в этом случае существует не более $C_{\rho\left(n\right)}^{3}$ вершин, удовлетворяющих данным условиям. 
\item \textbf{$supp\left(w\right) \not\subset X_{n}$}
\par Иными словами, носитель вершины $w$ не лежит в объединении носителей вершин первого и третьего типа. Возможны два случая: 
\begin{enumerate}
\item $\vert supp\left(w\right) \cap X_{n} \vert = 1.$ Поскольку вершина $w$ пересекается ровно с тремя вершинами независимого множества $I$, то данное равенство возможно только в том случае, если вершина $w$ пересекается с некоторым множеством вершин первого типа мощности 3. Пусть вершина $w$ имеет вид $\{x, y, z\}$, причем элемент $x$ принадлежит трём вершинам некоторого множества вершин первого типа. Тогда элемент $x$ можно выбрать не более, чем $n$ способами, элемент $y$ можно выбрать не более, чем $n$ способами, а элемент $z$ можно выбрать не более, чем тремя способами, так как иначе множество $I$ не было бы максимальным. Таким образом, в данном случае существует не больше $3 \cdot n^{2}$ вершин. 

\item $\vert supp\left(w\right) \cap X_{n} \vert = 2.$ Пусть вершина $w$ имеет вид $\{x, y, z\}$, причём элементы $x, y$ принадлежат множеству $X_{n}$. Пару элементов $x, y$ можно выбрать не более, чем $n^{2}$ способами, а элемент $z$ можно выбрать не более, чем тремя способами, так как иначе независимое множество $I$ не было бы максимальным. Таким образом, в данном случае существует не больше $3 \cdot n^{2}$ вершин. 
\end{enumerate}
Итого, число вершин, удовлетворяющих данному условию, не превосходит $6 n^{2}$.
\end{itemize}
\end{enumerate}

Таким образом, мы рассмотрели все возможные варианты взаимного расположения вершины $w$ и независимого множества. Просуммировав результаты всех пунктов мы получим, что существует не более чем $C_{\rho\left(n\right)}^{3} + 17n^{2} \leq \frac{\rho\left(n\right)^{3}}{6} + 20n^{2}$ вершин, пересекающихся ровно с тремя вершинами из независимого множества. Таким образом, лемма доказана.


\end{document}